\documentclass[journal]{IEEEtran}
\long\def\/*#1*/{}

\usepackage{amsmath,amsfonts}
\usepackage{lmodern}
\usepackage{graphicx,float}
\usepackage[switch]{lineno}
\usepackage{graphicx}

\title{\LARGE Search-Enhanced Instantaneous Frequency Detection Algorithm: A Preliminary Design}

\author{Phen Chiak See, and Marta Molinas}

\begin{document}

\maketitle


\begin{abstract}
This paper presents a method developed for finding sinusoidal components within a nonlinear non-stationary time-series data using Genetic Algorithm (GA) (a global optimization technique). It is called Search-Enhanced Instantaneous Frequency Detection (SEIFD) algorithm. The GA adaptively define the configuration of the components by simulating the solution finding process as a series of genetic evolutions. The start-time, end-time, frequency, and phase of each of these components are identified once convergence in the implementation is achieved. \end{abstract}


\section{Introduction}

\begin{figure}[h]
\centering
\includegraphics[scale=0.6]{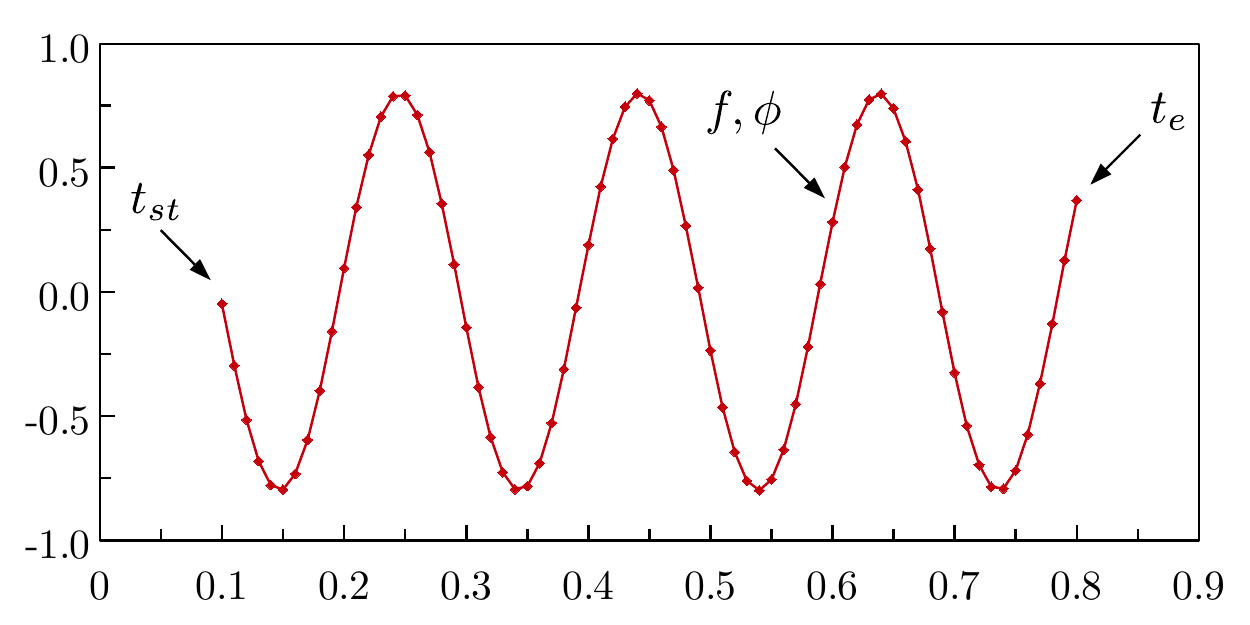}
\caption{A sinusoidal component within a time-series data.}
\label{fig:sinudoidal-component}
\end{figure}

The Search-Enhanced Instantaneous Frequency Detection (SEIFD) algorithm uses Genetic Algorithm (GA) \cite{Goldberg 1989, Whitley 1994} to identify the characteristics of sinusoidal components (see Fig. \ref{fig:sinudoidal-component}) within a nonlinear non-stationary time series data. The time-series data is modeled as a superposition of these components, each exist within a unique time range. It is represented as below.

\begin{equation}
f(t) = \sum_{n}^N a_n \sin(2\pi f_n t + 2 \phi_n t)
\end{equation}

\noindent where, $n$ is the index of each sinusoidal component; $a_n$ is its amplitude; $f_n$, is its frequency; and $\phi_n$ is its phase.

The parameters which define each of the components are encoded as a solution string consisting binary values. During implementation, the GA attempts to find their value by simulating it as a solution improvement process modeled as a series of genetic evolutions.

Because the problem is non-stationary, each of the sinusoidal components possesses unique start and end time (denoted respectively as $t_{st}^n$ and $t_e^n$). They are as well encoded in the solution string. In this way, the solver attempts to decompose the time-series data into its consituting sinusoidal components in order to reveal the instantaneous frequencies of the time-series data.

\section{The implementation}

\begin{figure}[htp]
\centering
\includegraphics[scale=0.6]{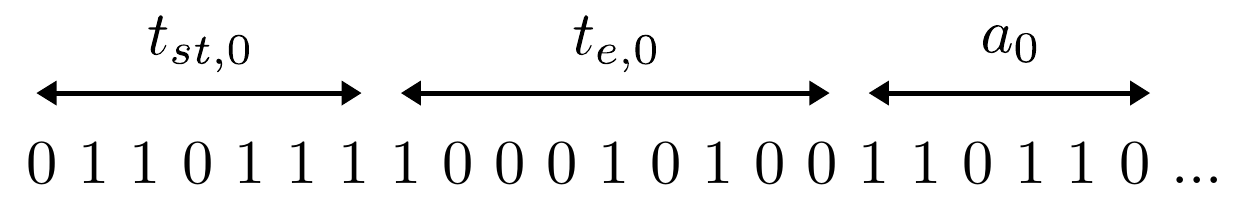}
\caption{A solution string in GA.}
\label{fig:sample-solution-string}
\end{figure}

The sinusoidal parameters are encoded as a string of binary values. The length of the binary string is determined based on the sum of the bit counts carried by each of the parameters computed with their respective lower and upper bound, as well as the step size of increment. The bit count for each of the parameters is computed as below.

\begin{equation}
\text{bit}(x) = \text{floor}\left[\text{ln} \left(\frac{x_{\text{ub}}-x_{\text{lb}}}{\Delta x}\right)\right] + 1
\end{equation}

\noindent where, bit($x$) is the bit count of parameter $x$; $x_{\text{ub}}$ and $x_{\text{lb}}$ respectively are the upper bound and lower bound of parameter $x$; and $\Delta x$ is the step size of each increment of parameter $x$ in the solution. The solution string is organized as shown in Fig. \ref{fig:sample-solution-string}.

The real value of each of the parameters can be decoded from the solution string using the following equation.

\begin{equation}
x = x_{\text{lb}} + \frac{x_{\text{ub}}-x_{\text{lb}}}{2^l-1}\left[\sum_{k=0}^{l-1} 2^k\text{bit}_{l-k}\right]
\end{equation}

\noindent where, $l$ is the length of bits in the solution string representing the value of $x$.

\begin{figure}[H]
\centering
\includegraphics[scale=0.6]{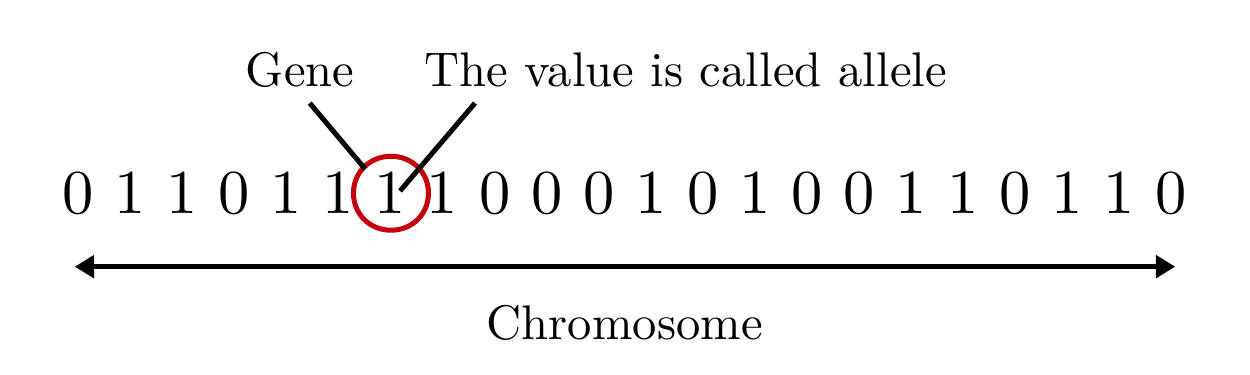}
\caption{The structure of chromosome in GA}
\label{fig:chromosome}
\end{figure}

The solution string is regarded as a chromosome (Fig. \ref{fig:chromosome}) in the simulated genetic evolution process. In the initial stage, GA generates a group (a population) of random chromosomes of varying quality computed with the following fitness function. 

\begin{equation}
F(x) = \text{min}\sum_{t=0}^T\left[ \sum_{n=0}^N a_n \sin(2\pi f_n t + 2 \phi_n t) - R(t)\right]
\end{equation}

\noindent where the $R(t)$ is the raw time-series data.

The fitness function describes the environment (the solution space) in which the chromosomes evolve. The implementation should achive continuous improvement (lower fitness function in every subsequent generation) until the fitness score converged at a targeted fitness score. The targeted fitness score of $F(x)$ is 0.

\begin{figure*}[t]
\centering
\includegraphics[scale=0.55]{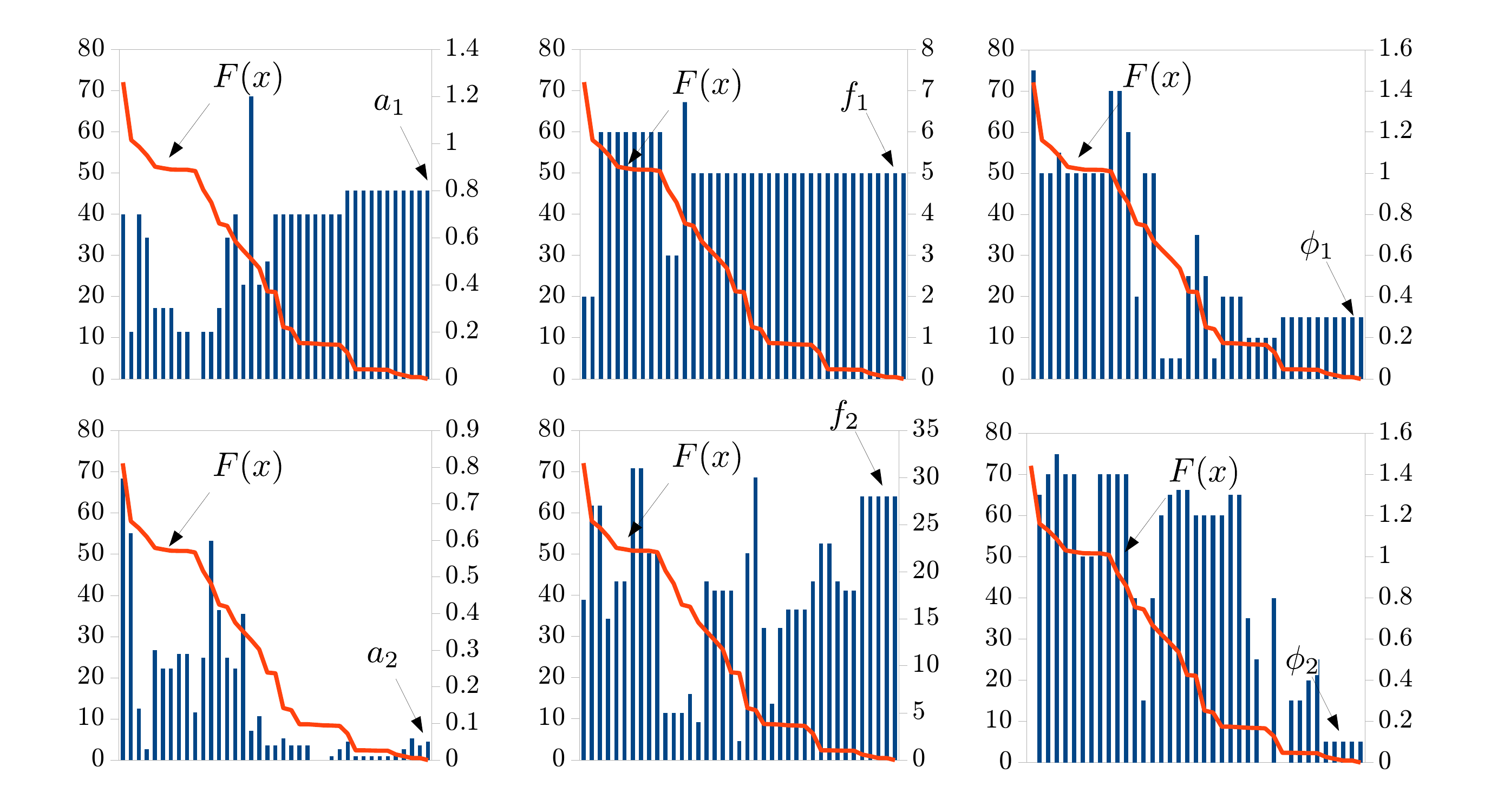}
\caption{The iteration-best value of $a_n$, $f_n$, and $\phi_n$ found by the GA.}
\label{fig:convergence}
\end{figure*}

Then, GA selects some chromosomes possessing good fitness score for creating new offsprings in the subsequent repetitions. In the process, the algorithm identify and retain good alleles that lead to the yielding of better fitness scores. Meanwhile, it modifies the inferior allele values using simulated biological crossover and mutation processes. The operation is repeated until the termination criterion is finally met.

Compare to Fourier transforms, the method adaptively identifies the start and end of such sinusoidal components (within the time horizon of the data), as well as their exact instantaneous attributes. Among others, its performance is not constrained by: i.) the need to ensure the completion of at least one oscillation cycle in all of the sinusoidal components; and ii.)  the existence of discontinuity within the data. The result fulfills the following generalized Fourier representation \cite{Norden Huang et al. 1998}:

\begin{equation}
X(t) = \sum_{j=1}^{n} a_j(t)\text{exp}\left( i\int \omega_i(t) \text{d}t \right)
\end{equation}

The identification of the quantity of sinusoidal components at every time step, $n_t$ can be solved as a Combinatorial Optimization (CO) problem (the GA is a well-known approach for solving CO problems (i.e., \cite{Anderson and Ferris 1994,Corcoran and Wainwright 1995})). To do so, the method could be implemented as below:

\begin{enumerate}

\item Establish a list of repeatable combination, $C$ of length $L$, which consist of a list of integer, $n_t \leq N^*$; where $L$ is the length of the data, $N^*$ is an integer representing the maximum quantity of sinusoidal components in the data.

\item Implement GA to determine the value of components in $C$, $c_i$.

\item Determine the $t_{st,j}$ and $t_{e,j}$ for each of the sinusoidal component based on $C$, and identify their configurations.

\item Improve $C$ until convergence is achieved.

\end{enumerate}

\section{Test case}

\begin{figure}[H]
\centering
\includegraphics[scale=0.55]{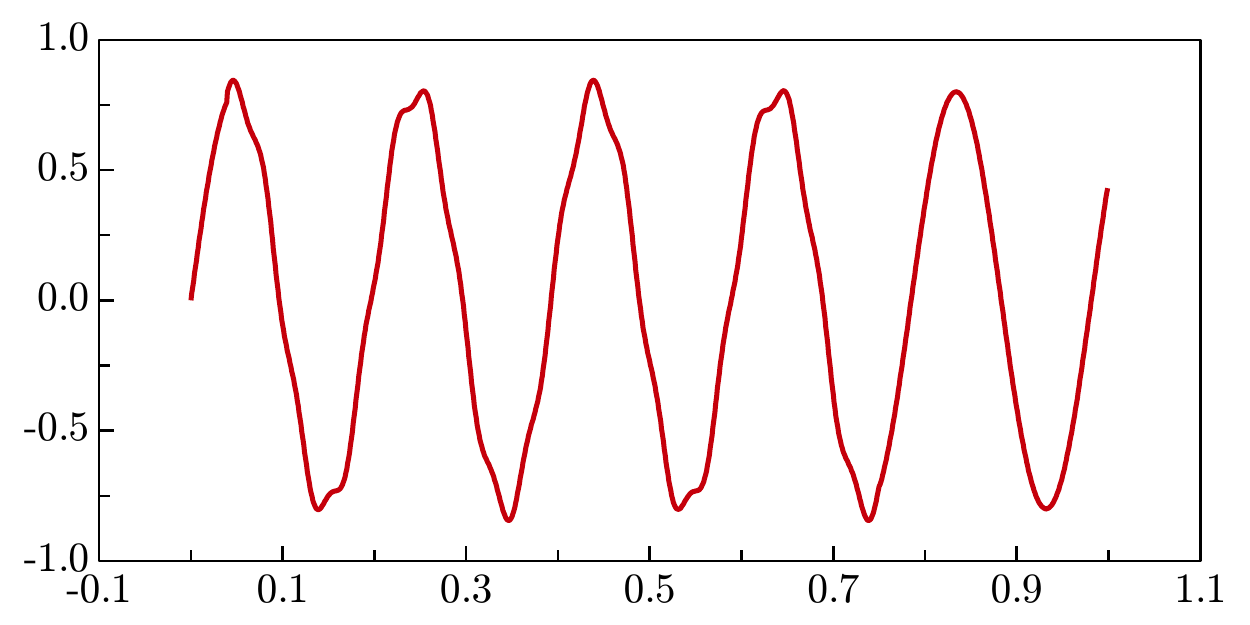}
\caption{A sample nonlinear non-stationary time-series data.}
\label{fig:sample-data}
\end{figure}

The method is tested on a time series data shown in Fig. \ref{fig:sample-data}. The algorithm is codified using the Python programming language.

\begin{figure}[H]
\centering
\includegraphics[scale=0.6]{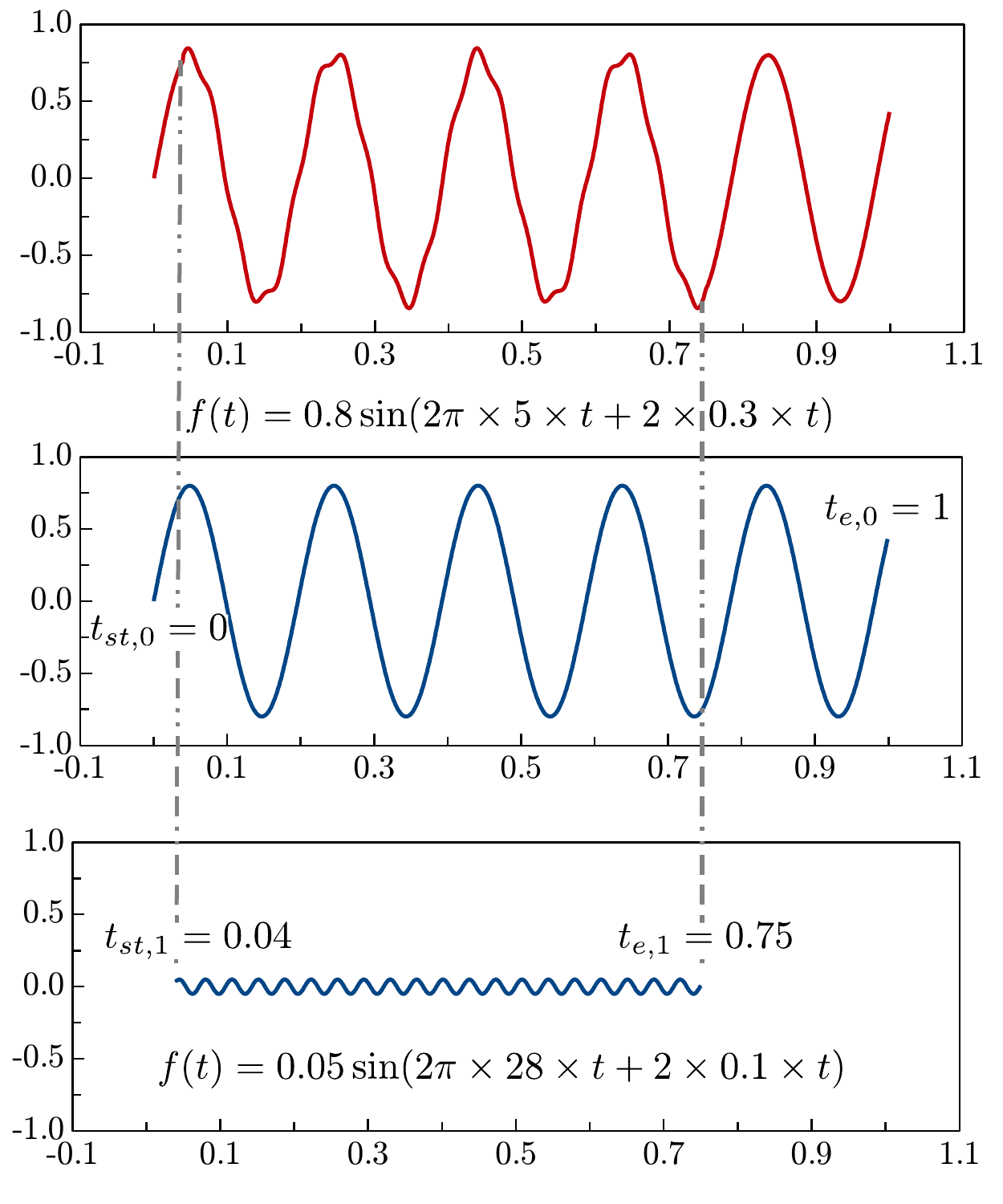}
\caption{The sinusoidal components in the test data.}
\label{fig:result}
\end{figure}

The iteration-best value of $a_n$, $f_n$, and $\phi_n$ found by the GA is shown in Fig. \ref{fig:convergence}. The changes in the parameter values which leads to the convergence of $F(x)$ is clearly seen in the figure. The final results obtained in the implementation is shown in Fig. \ref{fig:result}.

\end{document}